\documentstyle[hyperref,bm,cite,12pt]{article}
\RequirePackage{mathrsfs}
\hypersetup{colorlinks=true,pagebackref,pdftitle={Symmetries of vector exterior differential systems},pdfdisplaydoctitle=true}
\makeatletter
\def\ps@firstpage{\ps@plain
  \def\@oddhead{\article@logo\hss}
}
\def\article@logo{\vbox to\headheight{\@parboxrestore \small
\rightline{Journal of Nonlinear Mathematical Physics}
\rightline{Vol.~4 (1997) No.~1/2 (p.~89--97)}
\rightline{\url{http://staff.www.ltu.se/~norbert/home_journal/electronic/v4n1-2.html}}
   \par\vss
}}
\makeatother
\chardef\atcode=\catcode`\@
\catcode`\@=11
	\ifcase\@ptsize     % Fonts for 10pt style
\font\frak=eufm10   % Gothic for modules of differential forms and components of generator "x"
\font\sfrak=eufm7   % Subscript gothic to tensor product of modules
\font\bfrak=eufb10  % Infinitesimal operator "x" bold gothic
\font\sssm=cmss8    % Three-dimensional vectors in subscripts
\font\ssb=cmssbx10  % Three-dimensional vectors
\font\bb=msbm10     % Blackboard bold for real numbers field and matrices
	\or                 % Fonts for 11pt style
\font\frak=eufm10  scaled \magstephalf % Gothic for modules of differential forms and components of generator "x"
\font\sfrak=eufm7  scaled \magstephalf % Subscript gothic to tensor product of modules
\font\bfrak=eufb10 scaled \magstephalf % Infinitesimal operator x bold gothic
\font\sssm=cmss8   scaled \magstephalf % Three-dimensional vectors in subscripts
\font\ssb=cmssbx10 scaled \magstephalf % Three-dimensional vectors
\font\bb=msbm10    scaled \magstephalf % Blackboard bold for real numbers field and matrices
	\or                 % Fonts for 12pt style
\font\frak=eufm10  scaled \magstep1 % Gothic for modules of differential forms and components of generator "x"
\font\sfrak=eufm7  scaled \magstep1 % Subscript gothic to tensor product of modules
\font\bfrak=eufb10 scaled \magstep1 % Infinitesimal operator x bold gothic
\font\sssm=cmss8   scaled \magstep1 % Three-dimensional vectors in subscripts
\font\ssb=cmssbx10 scaled \magstep1 % Three-dimensional vectors
\font\bb=msbm10    scaled \magstep1 % Blackboard bold for real numbers field and matrices
	\fi
\catcode`\@=\atcode
\def\F{\hbox{\frak F}}
\def\I{\hbox{\frak I}}
\def\P{\hbox{\frak P}}
\def\S{\hbox{\frak S}}
\def\cF{{\mathscr F}}
\def\cI{{\mathscr I}}
\def\R{\hbox{\bb R}}
\def\A{\hbox{\bb A}}
\def\B{\hbox{\bb B}}
\def\x{\hbox{\bfrak x}}
\def\eisf#1{\hbox{\sssm#1}}
\def\bcdot{\mbox{\boldmath$\cdot$}}
\def\bpartial{\mbox{\boldmath$\partial$}}
\def\bmat#1{{\boldsymbol{#1}}}
\def\bkey#1{\hbox{\bf#1}}
\newcommand{\ra}{\rightarrow}
\newcommand{\df}{\mathrel{\stackrel{\rm def}{=}}}
\newcommand{\ul}{\underline}
\newcommand{\SS}{\scriptscriptstyle}
\newcommand{\Ss}{\scriptstyle}
\newcommand{\Hom}{\mathop{\rm Hom}}
\newcommand{\End}{\mathop{\rm End}}
\newcommand{\sol}{\mathop{\rm sol}}
\newcommand{\Far}[1]{\,$\F (B,#1)$}
\newcommand{\bw}{\bar{\wedge}}
\newcommand{\be}{\bar{\eta}}
\newtheorem{definition}{Definition}
\newtheorem{lemma}{Lemma}
\begin{document}
\title{\vspace*{1.5cm}Symmetries of vector exterior differential systems and
the inverse problem in second-order  Ostrohrads'kyj mechanics}
\author{R.~Ya.~Matsyuk\\\small email:~\texttt{romko.b.m@gmail.com}\\[-2\jot]\small \url{http://www.iapmm.lviv.ua/12/eng/files/st_files/matsyuk.htm}}
\date{\em Institute for Applied Problems in Mechanics and Mathematics}
\maketitle\thispagestyle{firstpage}
\begin{abstract}
Symmetries of variational problems are considered as symmetries of vector bundle valued exterior differential systems. This approach is then applied to third order ordinary variational equations of motion of the semi-classical spinning particle.
\end{abstract}
It is common in the similarity theories of differential equations
to consider transformations which  leave invariant some prolonged manifolds
in jet-spaces of appropriate order within which the integral manifolds
of a given system of partial differential equations lie~\cite{Ovs}.
From this invariance it follows that the solutions of the system in charge
transform but into some other solutions. However, we are not sure enough
that the prolonged manifold is ``densely'' covered by the integral manifolds.
If not, then some additional symmetries may occur. In this paper we
consider a symmetry of a differential equation be defined in general as
a generator of such a transformation, which carries every solution into
nothing more than some other solution.

Then we reformulate this  into slightly different language of the algebraic
invariance of some vector valued  exterior differential system. This has
an advantage that we can use Lagrange multiplies method to solve the problem.
In the case of ordinary differential equations the two approaches (algebraic
and the general one) coincide due to the completeness of Pfaff systems
(that is that every differential form, annulled by every solution of the
system expands by the Lagrange multiplies into the elements of the differential
ideal, generated by the system itself).

The Euler-Lagrange equations of the variational calculus (of arbitrary order)
naturally fall into the framework of vector bundle valued exterior differential systems
 when the dependent variables are globally segregated from the independent
ones~\cite{Kal}. This is the case of the field theory. In general such
segregation depends on the local chart
\,$\mu:\,M \ra\R^{p}\times\R^{q}$ \,
which constitutes a way of introducing a fibred structure in the manifold $M$.
In applying infinitesimal considerations we can use this structure every time the
manifold $M$ is not endowed by a fibred structure in an intrinsic manner.
If $Y$ denotes a fibred manifold over a base $Z$, we take
\,$Y_s\df J_sY$ \,
to mean the space of all the jets of order $s$ of cross-sections of the surmersion
\,$\pi:\,Y\ra Z$. Then
\,$\mu_{_{J}}:\,Y_s\ra J_s(\R^{p}\times\R^{q})\approx J_s(\R^p,\R^q)$ \,
will denote the standard $s$-order prolongation of the local chart $\mu$.
%%%%%%%%%%%%%%%%%%%%%%%%%%%%%%%%%%%%%
\paragraph*{1.~Vector bundle valued differential systems.}
Let
\,$E\ra B$,
\,$E_{\scriptscriptstyle1}\ra B_{\scriptscriptstyle1}$,
\,$F\ra X$,
\,$F_{\scriptscriptstyle1}\ra X_{\scriptscriptstyle1}$,
and
\,$E'\ra B$ \,
be some fibre bundles, let
\,$\Gamma (E)$ \,
denote the set of smooth cross-sections of the fibred manifold $E$ and
let
\,$\F (B)$ \,
stand for the ring of smooth functions defined on the manifold $B$. We shall
denote the tensor product of the inverse image of vector bundles with respect to
some manifold morphisms
\,$f:\,B\ra X$ \,
and
\,$f_{\scriptscriptstyle1}:\,B\ra X_{\scriptscriptstyle1}$ \,
by $\otimes_{_{B}} $. If the manifold $B$ is fibred over some manifold $N$,
the semi-basic differential forms over the manifold $B$ with values in the
vector bundle $F$ are defined as smooth cross-sections of the vector bundle
\,$\wedge T^{\ast}N\otimes_{_{B}}\F$.
They constitute an
\,$\F (B)$-submodule, denoted by
\,$\F_{_{B}}(N,F)$,
of the graded module
\,$\F(B,F)=\Gamma(\wedge T^{\ast}B\otimes_{_{B}}F)$.
Every bundle homomorphism
\,$h:\,F\ra F_{{\scriptscriptstyle1}}$ \,
over the base homomorphism
\,$\ul{h}:\,X\ra X_{\scriptscriptstyle1}$ \,
satisfying
\,$f_{\scriptscriptstyle1}=\ul{h}\circ f$ \,
obviously defines the moduli homomorphism
\,$h_{{\SS \#}}:\,\F(B,F)\ra\F(B,F_{\scriptscriptstyle1})$.
The canonical vector bundle pairing
\,$E\times\Hom(E,E')\ra E'$ \,
defines an exterior product
\,$\wedge:\,\F(B,E)\times\F(B,E^{\ast}\otimes E')\ra\F(B,E')$.
If we think of a cross-section
\,$\bmat{\varphi}\in\F^{\scriptscriptstyle0}(B,E^{\ast}\times E')$ \,
as a $B$-morphism
\,$\widehat{\phi}:\,E\ra E'$,
then
\,$\widehat{\phi}_{\SS \#}\bmat{\omega}=\bmat{\omega}\wedge\bmat{\varphi}$ \,
for any
\,$\bmat{\omega}\in\F(B,E)$.

The vector bundle
\,$\End(E)$ \,
is a bundle of algebras with the composition rule that of the superposition
of endomorphism. This fibrewise composition defines the structure of a
graded non-commutative algebra in the
\,$\F(B)$
module of
\,$\End (E)$-valued differential forms which we denote by
\,$\F(E,\End E)$.
The vector bundle $E$ may be viewed as the bundle of
\,$\End(E)$-moduli with respect to the pairing
\,$E\times\End(E)\ra E$.
This pairing makes the
$\F(B)$-module
\,$\F(B,E)$ \,
into a graded module over the algebra
\,$\F(B,\End E)$.

Now let
\,${\sf h}:\,E_{\scriptscriptstyle1}\ra E$ \,
be a vector bundle homomorphism over a base morphism
\,$g:\,B_{\scriptscriptstyle1}\ra B$ \,
and let
${\sf h}_{\scriptscriptstyle1}$
denote the corresponding induced $B_{\scriptscriptstyle1}$-homomorphism of vector bundles
\,${\sf h}_{\scriptscriptstyle1}:\,E_{\scriptscriptstyle1}\ra g^{-{\scriptscriptstyle1}}E$.
The dual vector bundle homomorphism
\,${\sf h}_{\scriptscriptstyle1}{}^\ast:\,(g^{-{\scriptscriptstyle1}}E)^\ast\ra E_{\scriptscriptstyle1}{}^\ast$ \,
denotes, by virtue of identification
\,$(g^{-{\scriptscriptstyle1}}E)^\ast\approx g^{-{\scriptscriptstyle1}}(E^\ast)$,
some $g$-comorphism
\,${\sf h}^\ast:\, g^{-{\scriptscriptstyle1}}(E^\ast)\ra E_{\scriptscriptstyle1}{}^\ast$.
If set
\,$E=TB$,
\,$E_{\scriptscriptstyle1}=TB_{\scriptscriptstyle1}$,
one gets the action
${\sf h}^{\SS \#}$ of ${\sf h}$
on the ordinary differential forms on~$B$. If
${\alpha}_g$ stand for the inverse image of the cross-section
\,$\alpha\in\F(B,\R)$,
then
\,${\sf h}^{\SS \#}{}{\alpha}=(\wedge {\sf h})^\ast\circ\alpha_g\in\F(B_{\scriptscriptstyle1},\R)$.
Given some cross-section
\,$\bmat{\omega}\in\F(B,E)$ \,
and recalling the notion of its inverse image
\,${\bmat{\omega}}_g\in\Gamma(\wedge T^\ast B\otimes_{B_{\scriptscriptstyle1}}E)$,
the differential form
\,$g^\star\bmat{\omega}\in\F(B_{\scriptscriptstyle1},E)$ \,
is being defined by the superposition of mappings:
\,$g^\star\bmat{\omega}=\big((\wedge Tg)^\ast\otimes {\sf id}\big)\circ\bmat{\omega}_g$.
If the morphism $g$ has a canonical prolongation to a vector bundle
$g$-comorphism
\,${\sf k}:\,g^{-{\scriptscriptstyle1}}E\ra E_{\scriptscriptstyle1}$,
then the notation
\,$g^{\SS \#}\bmat{\omega}$ \, will mean the differential form
\,${\sf k}_{\SS \#}g^\star\bmat{\omega}=\big((\wedge Tg)^\ast\otimes {\sf k}\big)\circ\bmat{\omega}$.

The functor $g^\star$ is consistent with the tensor product and contraction
operations. Let throughout this paragraph the manifold $B$ be parallelizable
and the module
\,$\Gamma (E)$ \,
be free. Then, if
\,$\bmat{\omega}=\alpha\otimes\bmat{\gamma}\in\F (B,E)\approx\F(B,\R)\otimes\Gamma (E)$ \,
and if
\,$\bmat{\Omega} =\beta\otimes\bmat{\varphi}\in\F (B,E^\ast\otimes E')\approx\F (B,\R)
\otimes \Gamma (E^\ast\otimes E')$,
then for the inverse images of the cross-sections
$\bmat{\omega}$,
$\bmat{\Omega}$, and
\,$\bmat{\omega}\wedge\bmat{\Omega} =(\alpha\wedge\beta )\otimes\langle{\bmat{\gamma}}_g,\bmat{\varphi}\rangle$ \,,
it is true that
\,$(\bmat{\omega}\wedge\bmat{\Omega})_g=({\alpha}_g\wedge {\beta}_g)\otimes
\langle {\bmat{\gamma}}_g,{\bmat{\varphi}}_g\rangle\in
\Gamma (\wedge T^\ast B{\otimes}_{_{B_1}}E')$.
Composing with
\,$(\wedge Tg)^\ast$ \,
one gets
\,$g^\star (\bmat{\omega}\wedge\bmat{\Omega} )=g^\star\bmat{\omega}\wedge g^\star\bmat{\Omega}$.

A one-form  $\bmat{\vartheta}\in\F^{\scriptscriptstyle1}(B,TB) $ may be thought of as a $B$-homomorphism of vector bundles
\,$\widehat{\theta}:\,TB\ra TB$; \,
then the dual homomorphism
\,$\widehat{\theta}^\ast:\,T^\ast B \ra TB$ \,
acts over the cross-sections from
\,$\Gamma (T^\ast B)\equiv\F^{\scriptscriptstyle1}(B,\R)$ \,
in an obvious manner and is being extended to
\,$\F(B,\R)$ \,
as a differentiation of degree $0$ acting trivially on the ring of functions.
 We denote this extended action by $\bmat{\vartheta}\bar{\wedge}$.
Again within local
considerations for a parallelizable $B$ and $\Gamma (E)$ being free,
we can identify
\,$\F(B,\End E)$ \,
with the tensor product
\,$\F (B,\R)\otimes_{_{\hbox{\sfrak F}\,{\Ss(B)}}}\Gamma (\End E)$ \,
and define
\,$\bmat{\vartheta}\bar{\wedge}(\alpha\otimes\bmat{\varphi} )=(\bmat{\vartheta}\bar{\wedge}\alpha)
\otimes{\sf id}(\bmat{\varphi})$.
Consider the ideal
\,$\I\equiv\I (B)\equiv\I (\bmat{\vartheta})=\sum\limits_{d>0}\I^{(d)}=
\bmat{\vartheta}\bar{\wedge}\F (B,\R)$ \,
in the algebra
\,$\F (B,\R)$.
Consider also the \,$\F (B,\R)$-submodule
\,$\I (B, \End E)=\bmat{\vartheta}\bar{\wedge}\F (B,\End E)$ \,
in the
\,$\F (B,\R)$-module
\,$\F (B,\End E)$.
We denote by
\,$\I (B,E)$ \,
the submodule
\,$\F (B,E)\wedge\I$ \,
in the module
\,$\F(B,E)$ \,
over the algebra
\,$\F (B,\R)$.
It is clear that
\,$\I (B,E)$ \,
is also a submodule over the algebra
\Far{\End E} \,
and that
\,$\I (B,E)=\F (B,E)\wedge\I(B,\End E)$.
In fact,
\,$\I(B,E)$ \,
is generated over
\Far{\R} \,
by the
\,$\F^{\scriptscriptstyle0}(B,\End E)$-submodule of one-forms
\begin{displaymath}
\P (B,E)=\bmat{\vartheta}\bar{\wedge}\F^{\scriptscriptstyle1}(B,E)=\F^{\scriptscriptstyle0}(B,E)\wedge\P,
\end{displaymath}
where
\,$\P\equiv\P(B)\equiv\I^{({\scriptscriptstyle1})}=\bmat{\vartheta}\bw\F^{\scriptscriptstyle1}(B,\R)$ \,
is an \,$\F(B)$-submodule, that is the ordinary Pfaff system.
%%%%%%%%%%%%%%%%%%%%%%%%%%%%%%%%%%%%%%%%%%%%%%%%%%%%%%%
\begin{definition}
An exterior differential system \,$\S\equiv\S(F)$ \,
(with values in the vector bundle $F$) is an
\Far{\End F}-submodule of the module
\Far{F}. A Pfaff system is an exterior differential system, generated by
one-forms. Let a  manifold~$Z$ be given. We call the germ of an immersion
\,$\sigma:\,Z\ra B$ \,
be a solution of the exterior differential system $\S$ if
\,$\sigma^\star\S=0$.
Two exterior differential systems (in general, over different manifolds) are
called equivalent if the sets of their solutions coincide (in general,
are isomorphic).
\end{definition}

Let over a manifold $B$ be given an exterior differential system $\S$
with values in the vector bundle $E$ and an exterior differential system
$\S'$ with values in the vector bundle $E'$ and let
\,$\sol (\S)$\,
denote the set of solutions of the system $\S$. If
\,$\S\wedge\F (B,E^\ast\otimes E')\supset\S'$,
then
\,$\sol (\S)\subset\sol (\S')$.
%%%%%%%%%%%%%%%%%%%%%%%%%%%%%%%%%%%%%%%%%%%%%%%%
\begin{definition}
Exterior differential systems $\S$ and $\S'$ are algebraically equivalent,
if at the same time both
\,$\S\wedge\F(B,E^\ast\otimes E')\supset\S'$
and
\,$\S'\wedge\F (B,E^{\prime\ast}\otimes E)\supset\S$.
\end{definition}
%%%%%%%%%%%%%%%%%%%%%%%%%%%%%%%%%%%%%%
\begin{lemma}
Two Pfaff systems which are algebraically equivalent, are equivalent.
\end{lemma}

This follows from the completeness of Pfaff systems.

Let $\widetilde{\S}$ denote the image of the exterior differential system $\S$
under the projection
\,$ j:\,\F (B,E)\ra\F (B,E)/\I (B,E)$.
In view of the inclusion
\,$\I (B,E)\wedge\F (B,E^\ast\otimes E')\subset\I (B,E')$ \,
the action $\wedge$ of the elements from
\Far{E^\ast\otimes E'} \,
over the quotient module is defined and if for this action \,$\S\wedge\F (B,E^\ast\otimes E')\supset \S'$,
then also
\,$\widetilde{\S}\wedge\F (B,E^\ast\otimes E')\supset\widetilde{\S'}$.
%%%%%%%%%%%%%%%%%%%%%%%%%%%%%%%%%%%%%%%%%%%%%%%%%%%%%%
\begin{definition}
A solution of the system
\,$(\S,\bmat{\vartheta})$ \,
is a germ of an immersion
\,$\sigma:\,Z\ra B$,
such that
\,$\sigma^\star j^{-{\scriptscriptstyle1}}(\widetilde{\S})=0$.
\end{definition}
%%%%%%%%%%%%%%%%%%%%%%%%%%%%%%%%%%%%%%%%%%%%%%%%%%%%

Let over the manifold $B_{\scriptscriptstyle1}$ be specified a  differential form
\,$\bmat{\vartheta}_{\scriptscriptstyle1}\in\F (B_{\scriptscriptstyle1},TB_{\scriptscriptstyle1})$ \,
such that
\,$g^\star\bmat{\vartheta} =\bmat{\vartheta}_{\scriptscriptstyle1}\bw\bmat{\varrho}\in\P (B_{\scriptscriptstyle1},g^{-{\scriptscriptstyle1}}TB)$ \,
for some
\,$\bmat{\varrho}\in\F (B_{\scriptscriptstyle1},TB)$.
Taking into account the definitions of the actions $g^\star$,
$\bmat{\vartheta}_{\scriptscriptstyle1}\bw$ and identifying the differential forms $\bmat{\vartheta}$,
$\bmat{\vartheta}_{\scriptscriptstyle1}$, and $\bmat{\varrho}$ with the corresponding vector bundle
homomorphisms $\widehat{\theta}$, $\widehat{\theta}_{\scriptscriptstyle1}$ and
\,$\widehat{\rho}:\, TB_{\scriptscriptstyle1}\ra g^\ast TB$ \,
we have
\,$\widehat{\theta}_g\circ (Tg)_{\scriptscriptstyle1}=\widehat{\rho}\circ\widehat{\theta}_{\scriptscriptstyle1}$,
where $\widehat{\theta}_g$ is the inverse image of the homomorphism $\widehat{\theta}$
with respect to the morphism $g$. For the conjugated homomorphisms we
also have
\,$(Tg)^\ast\circ\widehat{\theta}_g{}^\ast =\widehat{\theta}_{\scriptscriptstyle1}{}^\ast\circ\widehat{\rho}^{\,\ast}$.

In view of the definition of the action $\widehat{\rho}^{\,\SS \#}$
on the differential form
\,$\alpha\in\F^{{\scriptscriptstyle1}}(B,\R)$ \,
it holds that
\begin{displaymath}
g^\star (\bmat{\vartheta}\bw\alpha )=(Tg)^\ast\circ\widehat{\theta}_g{}^\ast\circ\alpha_g=
\widehat{\theta}_{\scriptscriptstyle1}{}^\ast\circ\widehat{\rho}^{\,\ast}\circ\alpha_g=
\bmat{\vartheta}_{\scriptscriptstyle1}\bw
\widehat{\rho}^{\,\ast}\alpha\in\P_{\scriptscriptstyle1}\equiv\P (B_{\scriptscriptstyle1}).
\end{displaymath}
It follows that
\,$g^\star\P\subset\P_{\scriptscriptstyle1}$,
thus
\,$g^\star\I\subset\I_{\scriptscriptstyle1}$ \,
and also
\begin{displaymath}
g^\star\I (B,E)=g^\star\F (B,E)\wedge g^\star\I\subset\F (B_{\scriptscriptstyle1},g^{-{\scriptscriptstyle1}}E)
\wedge\I (B_{\scriptscriptstyle1},g^{-{\scriptscriptstyle1}}E),
\end{displaymath}
so the action $g^\star$ of the homomorphism $g$ over the quotient module
\Far{E/\I (B,E)} \,
is defined.
%%%%%%%%%%%%%%%%%%%%%%%%%%%%%%%%%%%%%%%%%%%%%%%%%%%%
\paragraph*{2.~The Lie derivative.}
Let
\,$w:\,B\ra W$\,
be a morphism of manifolds and let $w_{t}$ denote its deformation. Then
\,$\dot{w}_t(0):\,B\ra TW$ \,
is a lift of the morphism $w$. If the manifold $W$ is fibred over the
manifold $B$ and if every $w_{t}$ is a cross-section, then
\,$\dot{w}_{t}(0)$ \,
belongs to the vertical tangent  bundle
\,$VW\subset TW$.
Let $\xi$ and $\bar{\eta}$ be some vector fields on the manifolds $B$ and $W$
respectively. The Lie derivative~\cite{Kol} of the morphism $w$
with respect to the pair of vector fields $\xi$, $\bar{\eta}$ is defined by
 the expression
\,${\bf L}_{\xi,\be}w =(Tw)\circ\xi-\be\circ w$.
If $w$ is a cross-section and if the vector field $\be$ is projectible onto the
vector field $\xi$, then the Lie derivative
\,${\bf L}_{\xi,\be}w$ \,
belongs to $VW$. Let under this assertions
\,$\exp_{_{B}}t\be$ \,
denote the one-parameter induced family of diffeomorphism over $B$,
corresponding to
\,$\exp t\be$ \,
and let
\begin{displaymath}
w_t=\exp_{_{B}}(-t\be)\circ w_{\exp t\xi}=\exp(-t\be)\circ w\circ\exp t\xi.
\end{displaymath}
Differentiating at
\,$t=0$ \,
we obtain
\begin{displaymath}
\dot{w}_t(0)=-(\exp t\be)\circ w +(Tw)\circ(\exp t\xi)\dot{{}}=
{\bf L}_{\xi,\be}\,w.
\end{displaymath}
If the fibre bundle $W$ is a vector bundle, then
\,$VW\approx W\times_{_{B}}W$ \,,
the first component of
\,$\dot{w}_t(0)$ \,
coincides with the initial cross-section $w$, and the second one is identified with the
Lie derivative
\,${\bf L}_{\xi,\be}\,w$ \,
itself. In a special case setting
\,$W=\wedge T^\ast B\otimes E$ \,
and taking $\be$ to be build up from the standard lift to the cotangent bundle of the
vector field $\xi$ together with some vector field $\eta$ on $E$,
projectible into $\xi$, we agree with truncated notation of
\,${\bf L}_{\xi,\eta}\df{\bf L}_{\xi,\be}$.
%%%%%%%%%%%%%%%%%%%%%%%%%%%%%%%%%%%%%%%%%%%%%%%%%%%%%%%%
\begin{definition}
A diffeomorphism
\,$g:\,B\ra B$ \,
is called a symmetry ({\it vs}  an algebraic symmetry) of the exterior differential
system $\S(E)$, if the exterior differential systems
$g^\star\S$ and $\S$ itself are equivalent ({\it vs} algebraically equivalent).
A vector field $\xi$ on the manifold $B$ is called an infinitesimal
symmetry ({\it vs} an algebraic infinitesimal symmetry) of the exterior differential
system $\S$ if for all $t$ the transformation
\,$\exp t\xi$ \,
 is a symmetry ({\it vs} an algebraic symmetry) of the system $\S$.
\end{definition}
%%%%%%%%%%%%%%%%%%%%%%%%%%%%%%%%%%%%%%%%%%%%%%%%

That $\xi$ is an algebraic infinitesimal symmetry of the exterior differential
system $\S$ means that
\begin{displaymath}
(\exp t\xi)^\star\S\subset\S\wedge\F\big(B,E^\ast\otimes(\exp t\xi)^{\SS-1}E\big).
\end{displaymath}
In particular, whatever the vector field
\,$\eta\in\Gamma(TE)$,
projectible into the vector field $\xi$ be, the inclusions
{
\setlength{\arraycolsep}{0pt}
\begin{eqnarray*}
(\exp -t\eta)_{_{B^{\SS\#}}}(\exp t\xi)^{\SS-1}\S&\subset&
\S\wedge\F\big(B,E^\ast\otimes
(\exp t\xi)^{\SS-1}E\wedge\F^{\scriptscriptstyle0}(B,{(\exp t\xi)^{\SS-1}}^\star(E^\ast\otimes E)\big)\\
&\subset&\S\wedge\F(B,E^\ast\otimes E)=\S
\end{eqnarray*}
}
hold.
By differentiating with respect to the parameter $t$ we conclude that if
$\xi$ is the infinitesimal symmetry of the system $\S$, then
\begin{displaymath}
{\bf L}_{\xi,\eta}\S=\frac{d}{dt}(\exp t\xi)^\star\S(0)\subset\S
\end{displaymath}
for every vector field $\eta$, projectible into the field $\xi$.
In the case of a trivial vector bundle $E$, one can put
\,$\eta=(\xi,0)$ \,
and use the notation
\,${\bf L}_{\xi,0}$ \, in place of
\,${\bf L}_{\xi,(\xi,0)}$.
%%%%%%%%%%%%%%%%%%%%%%%%%%%%%%%%%%%%
\begin{definition}
A diffeomorphism
\,$g:\,B\ra B$ \,
is called a symmetry ({\it vs} algebraic symmetry) of the system
\,$(\S,\bmat{\vartheta})$,
if
\,$g^\star\I\subset\I$ \,
and if the exterior differential systems
\,$g^\star\widetilde{\S}$ \,
and $\widetilde{\S}$ itself are equivalent ({\it vs} algebraically equivalent).
\end{definition}
%%%%%%%%%%%%%%%%%%%%%%%%%%%%%%%%%%%%%%%%%%%

An (algebraic symmetry of system
\,$(\S,\bmat{\vartheta})$ \,
is symmetry ({\it vs} an algebraic symmetry) of the exterior differential system
\,$\S+\I(B,E)$ \,
and {\it vice versa}. In that case
\,${\bf L}_{\xi}\I\subset\I$ \,
and
\,${\bf L}_{\xi,\eta}\widetilde{\S}\subset\widetilde{\S}$, that is
\begin{equation}\label{1}
{\bf L}_{\xi,\eta}\S\subset\S+\I(B,E)\,.
\end{equation}
%%%%%%%%%%%%%%%%%%%%%%%%%%%%%%%%%%%%%%%%%%
\paragraph*{3.~Symmetries of the Euler-Lagrange equations.}
On the manifold $Y_{s+{\scriptscriptstyle1}}$ of the jets of order $s+1$ of the cross-sections
of the fibred manifold
\,$\pi:\,Y\ra Z$ \,
there exists a canonical contact differential form
\,$\bmat{\vartheta}_s\in\F_{s+{\scriptscriptstyle1}}(Y_s,V_s)$, semi-basic with respect to the
projection
\,$\pi_s^{s+{\scriptscriptstyle1}}:\,Y_{s+{\scriptscriptstyle1}}\ra Y_s$ \,
which takes values in the vector bundle $V_s$ of vertical tangent vectors
to the surmersion
\,$\pi^{s}:\,Y_s\ra Z$ \,
and such that Pfaff system
\,$\P_s=\bmat{\vartheta}_s\bw{\F}_{s+{\scriptscriptstyle1}}^{\scriptscriptstyle1}(Y_s,\R)$ \,
is nothing but the Cartan co-distribution.
(We use subscript $s$ instead of $Y_s$ where possible).
For an open set
\,$U_s\subset Y_s$ \,
and applying the projection
\,$\pi_v^s:\,Y_s\ra Y_v$ \,
let us put
\,$\I_v(U_s)=\I(\left.\pi_v^s\right|_{_{U_s}}{\!\!\!\!}^\star\,\bmat{\vartheta}_v)$ \,
and let
\,$\upsilon_{_{J}}:\,Y_s\ra J(\R^p,\R^q)$ \,
be the canonical lift of a local fibred chart $\upsilon$ on the manifold $Y$.
A lagrangian is an element $\widetilde{\lambda}$ of the quotient sheaf, generated by
the quotient moduli
\,$\F_{_{U_r}}^p(Y,\R)/\I_{\scriptscriptstyle0}^{(p)}(U_r)$,
which we shall briefly denote by
\,${\cF}_{\,r}^{\,p}(Y,\R)/{\cI}_{\,\scriptscriptstyle0}^{\,(p)}(Y_r)$.
(We denote  the corresponding sheaves by the calligraphic characters).
There can always be found, at least locally, a  representative $\lambda$, semi-basic with
respect to the projection
\,$\pi^r:\,Y_r\ra Z$,
that is
\,$\lambda\in{\cF}_r{}^p(Z,\R)$.
The Euler-Lagrange expressions which correspond to the Lagrange density
\,$\upsilon_{_J}{}^{-{\scriptscriptstyle1}}{}^\star\lambda\in{\cF}_r{}^p(\R^p,\R)$ \,
naturally take the shape of the components of the local expression of
some differential $p$-form
\begin{equation}\label{2}
\bmat{\varepsilon}\in{\cF}_{\,2r}^{\,p}(Z,V_{\scriptscriptstyle0}{}^\ast)
\end{equation}
corresponding to the lagrangian $\lambda$~\cite{Kaz,MatStud}.

{\em
The symmetries of the Euler-Lagrange equations are nothing else but the
symmetries of the system
\,$(\S_{_{\bmat{\varepsilon}}},\bmat{\vartheta}_{2r-{\scriptscriptstyle1}})$,
where the module $\S_{_{\bmat{\varepsilon}}}$ is generated by the vector bundle valued differential
form $\bmat{\varepsilon}$.
}
%%%%%%%%%%%%%%%%%%%%%%%%%%%%%%%%%%%%%%%%%%%%%%%%%%%%%%%%%%%%%%%%%%%%%%%%%%%%
From here on we shall deal with variational calculus in one independent
variable. As an example we consider a system of third-order ordinary
differential equations
\begin{equation}\label{18}
{{\sf E}}_{a}=0\,.
\end{equation}
Let \,$Y=\R^p\times\R^q$, and let the
canonical coordinates in the manifold $J_{r}(\R^{1};\R^{q})$ be denoted by
$t, \ \hbox{\ssb x}=({\sf x}^{a}), \ \hbox{\ssb v}=({\sf v}^{a}),
\ \hbox{\ssb v}^{{\prime}}=({{\sf v}'}^{a}), \, \dots, \hbox{\ssb v}^{\scriptscriptstyle(r-1)}
=({{\sf v}^{\scriptscriptstyle(r-1)}}^{a})$.

This corresponds to \,$p=1$, $k=q$ \, in~(\ref{2}).
We introduce a vector valued differental one-form
\begin{equation}\label{19}
{\bmat{\epsilon}}={{\sf E}}_{a}\,{\bf d}{{\sf x}}^{a}\otimes{\bf d}t\,,
\end{equation}
where the expressions ${{\sf E}}_{a}$ are the Euler-Poisson
expressions. Applying the general criterion of \cite{Tul} for an arbitrary system
of differential equations to be a system of Euler-Poisson equations, it was
established in \cite{Thesis,MatMet20} that the vector expression ${\hbox{\ssb E}}=\{{{\sf E}}_{a}\}$
in (\ref{18}) must take the shape of
\begin{equation}\label{20}
{\hbox{\ssb E}}=\A{\,\bkey.\,}{\hbox{\ssb v}}^{{{\prime}}{{\prime}}}{\,+\,}({\hbox{\ssb v}}^{{\prime}}{\!\bkey.\,}{\bpartial}_{\eisf v})\,
\A{\,\bkey.\,}{\hbox{\ssb v}}^{{\prime}}{\,+\,}\B{\,\bkey.\,}{\hbox{\ssb v}}^{{\prime}}{\,+\,}{\hbox{\ssb c}}\,,
\end{equation}
where  the skew-symmetric matrix $\A$, the matrix $\B$, and the column vector
${\hbox{\ssb c}}$ depend on the variables $t$, ${\hbox{\ssb x}}$, ${\hbox{\ssb v}}={d{\hbox{\ssb x}}}/dt$,
and satisfy the following system of partial differential equations in $t$, ${{\sf x}}^a$,
and ${{\sf v}}^a$~\cite{Thesis,Dokl}
{\openup1\jot
\begin{equation}\label{21}
\begin{array}{c}
       \partial_{_{_{_{{\hbox{\sssm v}}}}}}{\!}_{[a}{}{{\sf A}}_{bc]}=0 \\
       2\,{{\sf B}}_{[ab]}-3\,{\bf D_{_{\bmat 1}}}{\kern.01667em}{{\sf A}}_{ab}=0 \\
       2\,\partial_{_{_{_{{\hbox{\sssm v}}}}}}{\!}_{[a}{}{{\sf B}}_{b]c}
-4\,\partial_{_{_{_{{\hbox{\sssm x}}}}}}{\!}_{[a}{}{{\sf A}}_{b]c}
+{\partial_{_{_{_{{\hbox{\sssm x}}}}}}{\!}_{c}}{\,}{{\sf A}}_{ab}
+2\,{\bf D_{_{\bmat 1}}}{\kern.01667em}{\partial_{_{_{_{{\hbox{\sssm v}}}}}}{\!}_{c}}{\,}{{\sf A}}_{ab}=0 \\
       {\partial_{_{_{_{{\hbox{\sssm v}}}}}}{\!}_{(a}}{}{{\sf c}}_{b)}
-{\bf D_{_{\bmat 1}}}{\kern.01667em}{{\sf B}}_{(ab)}=0\\
%%%%%%%%%%%%%%%
       2\,{\partial_{_{_{_{{\hbox{\sssm v}}}}}}{\!}_{c}}{\,}\partial_{_{_{_{{\hbox{\sssm v}}}}}}{\!}_{[a}{}{{\sf c}}_{b]}
-4\,\partial_{_{_{_{{\hbox{\sssm x}}}}}}{\!}_{[a}{}{{\sf B}}_{b]c}
+{{\bf D_{_{\bmat 1}}}}^{2}{\,}{\partial_{_{_{_{{\hbox{\sssm v}}}}}}{\!}_{c}}{\,}{{\sf A}}_{ab}
+6\,{\bf D_{_{\bmat 1}}}{\kern.0334em}\partial_{_{_{_{{\hbox{\sssm x}}}}}}{\!}_{[a}{}{{\sf A}}_{bc]}=0 \\
       4\,\partial_{_{_{_{{\hbox{\sssm x}}}}}}{\!}_{[a}{}{{\sf c}}_{b]}
-2\,{\bf D_{_{\bmat 1}}}{\kern.0334em}\partial_{_{_{_{{\hbox{\sssm v}}}}}}{\!}_{[a}{}{{\sf c}}_{b]}
-{{\bf D_{_{\bmat 1}}}}^{3}{\,}{{\sf A}}_{ab}=0\,.
\end{array}
\end{equation}
}
In~(\ref{21}) $\mathbf{D_{{s}}}$ denotes the generators of Cartan distribution,
{\setlength{\arraycolsep}{0pt}
\begin{eqnarray*}
{\bf D_{_{\bmat 2}}}={\hbox{\ssb v}}^{{\prime}}{\!\bkey.\,}{\bpartial}_{\eisf v}
{\,+\,}&{\bf D_{_{\bmat 1}}}&\,,\\
&{\bf D_{_{\bmat 1}}}&=\partial_{t}{\,+\,}{\hbox{\ssb v}}{\,\bkey.\,}{\bpartial}_{\eisf x}\,.
\end{eqnarray*}

Let $\bmat{\vartheta}_{_{\bmat 2}}$, $\bmat{\vartheta}_{_{\bmat 3}}$ denote the canonical contact
forms
\begin{eqnarray*}
\bmat{\vartheta}_{_{\bmat 3}}
=\frac\partial{\partial {{{\sf v}}'}^{a}}\otimes({\bf d}{{\sf v}'}^{a}
-{{\sf v}''}^{a}{\bf d}t){\,+\,}&\bmat{\vartheta}_{_{\bmat 2}}&\,,  \\
&\bmat{\vartheta}_{_{\bmat 2}}&
=\frac\partial{\partial {{\sf v}}^{a}}\otimes({\bf d}{\sf v}^{a}
-{{\sf v}'}^{a}{\bf d}t)
{\,+\,}\frac\partial{\partial {{\sf x}}^{a}}\otimes({\bf d}{\sf x}^{a}
-{\sf v}^{a}{\bf d}t) \,.
\end{eqnarray*}
Along with the differential form ${\bmat{\epsilon}}$,
we introduce another one,
\begin{eqnarray*}
{\underline{\bmat{\epsilon}}}
={\sf A}_{ab}\,{\bf d}{\sf x}^{a}\otimes{\bf d}{{\sf v}'}^{b}
{\,+\,}&{\sf K}&_{a}\,{\bf d}{\sf x}^{a}\otimes{\bf d}t\,,\\
&\hbox{\ssb K}&
=({\hbox{\ssb v}}^{{\prime}}{\!\bkey.\,}{\bpartial}_{\eisf v})\,
\A{\,\bkey.\,}{\hbox{\ssb v}}^{{\prime}}{\,+\,}\B{\,\bkey.\,}{\hbox{\ssb v}}^{{\prime}}{\,+\,}{\hbox{\ssb c}}\,.
\end{eqnarray*}
}
Exterior differential systems, generated by the forms ${\bmat{\epsilon}}$ and
${\underline{\bmat{\epsilon}}}$, are equivalent:
\begin{equation}\label{33}
{\underline{\bmat{\epsilon}}}-{\bmat{\epsilon}}
=
\bmat{\vartheta}_{_{\bmat 3}}\bw({\sf A}_{ab}\,{\bf d}{\sf x}^{a}\otimes{\bf d}{{\sf v}'}^{b})\,.
\end{equation}

Now it is time to put in the concept of {\it symmetry}. Let
\begin{equation}\label{22}
{\x}=\hbox{\frak t}\frac\partial{\partial t}{\,+\,}{\hbox{\frak x}^{a}}\frac\partial{\partial{\sf x}^{a}}
\end{equation}
denote a generator of some local group of transformations of the manifold \,$\R\times\R^q$,
its successive prolongations to the space $J_{s}(\R;\R^q)$
denoted by ${\x}_{_{\bmat s}}$,
\[
{\x}_{_{\bmat 2}}=\hbox{\frak v}^{a}\frac\partial{\partial{{\sf v}'}^{a}}
{\,+\,}{\x}_{_{\bmat 1}}\,.
\]

The demand that the exterior differential system, generated by the vector
valued differential form ${\underline{\bmat{\epsilon}}}$, be invariant under the
infinitesimal transformation ${\x}$ incarnates in consistency with
(\ref{1}) into the following
equation
\begin{equation}\label{23}
{\bf L}({\x}_{_{\bmat 2}})({\underline{\bmat{\epsilon}}})={\bf\Xi}{\,\bkey.\,}
{\underline{\bmat{\epsilon}}}{\,+\,}\bmat{\vartheta}_{_{\bmat 2}}\bw{\bmat{\omega}}\,,
\end{equation}
where the elements of a matrix
${\bf\Xi} \in \F^{0}\big(J_{2}(\R;\R^{q});\Hom({\R^{q}}^{\ast};{\R^{q}}^{\ast})\big)$
and the coefficients of a $\pi^{2}_{1}$-horizontal ${\R^{q}}^{\ast}$-valued one-form
${\bmat{\omega}}\in \F^{1}_{2}\big(J_{1}(\R;\R^{q});{\R^{q}}^{\ast}\big)$
depend upon the variables $t$, $\hbox{\ssb x}$, $\hbox{\ssb v}$, and ${\hbox{\ssb v}}^{{\prime}}$.
Both ${\bf\Xi}$ and ${\bmat{\omega}}$ play the role of Lagrange multipliers.
Splitting equation~(\ref{23})
with respect to independent differentials
${\bf d}t$, ${\bf d}{\sf x}^{a}$, ${\bf d}{\sf v}^{a}$, and
${\bf d}{{\sf v}'}^{a}$, gives the following system of partial differential
equations
{\openup1\jot
\begin{eqnarray}\label{24}
{\bf L}({\x}_{_{\bmat 1}}){\sf A}_{ab}&=&{\Xi_{a}}^{c}\,{\sf A}_{cb}
-{\sf A}_{ac}{\displaystyle\frac\partial{\partial{{\sf v}'}^{b}}\hbox{\frak v}^{c}}\\
{\bf L}({\x}_{_{\bmat 2}}){\sf K}_{a}&=&{\Xi_{a}}^{b}\,{\sf K}_{b}
-{\sf A}_{ab}{\bf D_{_{\bmat 2}}}\hbox{\frak v}^{b}
-{\sf K}_{a}{\bf D_{_{\bmat 1}}}\hbox{\frak t}\,.
\end{eqnarray}
}
\subparagraph*{Variational problems in parametric form.}
A variational problem in parametric form is a variational problem, posed on the
manifold \,$J_r(Z,M)$, so we have to put \,$Y=Z\times M$, the dimension of the
manifold $M$ equal to \,$p+q$ with $p=1$, and \,$k=p+q=1+q$ \, in~(\ref{2}).
Canonical coordinates in the manifold $J_{r}(\R;\R^{1+q})$ are denoted by
$\tau, \ \bmat x=(x^{\rho}), \ \bmat u=(u^{\rho}), \ \bmat{\dot{ u}}
=({\dot{u}}^{\rho}), \, \dots,{\stackrel{\scriptscriptstyle r-1}{\bmat u}}
=({\stackrel{\scriptscriptstyle r-1}{u}}^{\rho})$.
The manifold of $r^{\rm th}$-order velocities, $T^{1}_{r}M$, is defined as
$T^{1}_{r}M=J_{r}(\R;M)_{0}$. There exists an obvious isomorphism
$J_{r}(\R;M)\approx\R\times T^{1}_{r}M$.
Coordinates in the manifold $T^{1}_{r}M$ are denoted by
$\bmat x, \bmat u,
\bmat{\dot{ u}}, \dots,{\stackrel{\scriptscriptstyle r-1}{\bmat u}}$.
If $M$ is a (pseu\-do-\negthinspace\nolinebreak) \negthinspace\negthinspace Eu\-c\-li\-d\-ean $n$-dimensional space
(of an arbitrary signature), the Hodge operator ``$\ast$'' is defined as
$(\ast w)_{\rho_{k+1} \dots \rho_{n}}=\frac1{(n-k)!}
e_{\rho_{1} \dots \rho_{n}}w^{\rho_{1} \dots \rho_{k}}$.

Consider now a variational problem in parametric form, set by a Lagrangian
$$
\ell(\tau, x^{\rho}, u^{\rho}, \dots,{\stackrel{\scriptscriptstyle r-1}{u}}{}^{\rho}){\bf d}\tau
$$
on the space $J_{r}(\R;M)$. As long as we constrain ourselves only
to the case of autonomous Euler-Poisson equations,
\begin{equation}\label{25}
{\cal E}_{\rho}=0\,,
\end{equation}
the differential form
\begin{equation}\label{27}
\bmat{\varepsilon}={\cal E}_{\rho}\,{\bf d}x^{\rho}\otimes{\bf d}\tau\,,
\end{equation}
may {\em globally\/} be deprived of the factor ${\bf d}\tau$, constituting
thus a {\em globally\/} defined $T^{\ast}M$-valued density
\begin{equation}\label{28}
\bmat e={\cal E}_{\rho}{\bf d}x^{\rho}\,.
\end{equation}

Let $C^1_rM$ denote the manifold of the $r^{\rm th}$ order contact germs
of one-dimensional submanifolds of $M$.
The projection $\wp:T^{1}_{r}M\setminus0\to C^{1}_{r}M$ can be employed to
generate an autonomous variational problem set over $T^{1}_{r}M$ from every one
variational problem over $C^{1}_{r}M$.

\begin{lemma}[\cite{Kal}]
In terms of a local chart, if in (\ref{19}) the local semi-basic differential
form ${\bmat{\epsilon}}$ corresponds to the Lagrangian
$$
{\lambda}=L{\bf d}t\,,
$$
then the vector valued density
\begin{equation}\label{29}
\bmat e=-u^{a}({\sf E}_{a}\circ\wp){\bf d}x^{\scriptscriptstyle0}
{\,+\,}u^{\scriptscriptstyle0}({\sf E}_{a}\circ\wp){\bf d}x^{a}
\end{equation}
corresponds to the Lagrangian
$$
\ell(\tau, x^{\rho}, u^{\rho}, \dots,{\stackrel{\scriptscriptstyle r-1}{u}}{}^{\rho}){\bf d}\tau
={\cal L}(x^{\rho}, u^{\rho}, \dots,{\stackrel{\scriptscriptstyle r-1}{u}}{}^{\rho}){\bf d}\tau
$$
with the Lagrange function
\begin{equation}\label{30}
{\cal L}=u^{\scriptscriptstyle0}L{\,\circ\,}\wp\,.
\end{equation}
\end{lemma}

\paragraph{4.~The classical spinning particle.}
Let in (\ref{22}) generator ${\x}$ correspond to the (pseu\-do-\negthinspace\nolinebreak) \negthinspace\negthinspace or\-tho\-go\-nal
transformations of a four-dimensional (pseu\-do-\negthinspace\nolinebreak) \negthinspace\negthinspace Eu\-c\-li\-d\-ean space. In this case,
equations (\ref{24}) have no solutions of the third order, that is to say,
no invariant system of Euler-Poisson equations (\ref{18}) consists of at least one equation
of the third order.

Nevertheless, if we allow a vector parameter $\bmat s=(s^{\scriptscriptstyle0},\hbox{\ssb s})$
(that is, a constant
quantity, transforming as a four-vector under the action of the (pseu\-do-\negthinspace\nolinebreak) \negthinspace\negthinspace or\-tho\-go\-n\-al
group) enter into expressions (\ref{20}), then equations (\,\ref{21}\,\&\ref{24}\,)  turn out to
possess a family of solutions, which depend on $\bmat s$ and actually contain some
third-order derivatives. In~\cite{Thesis,Dokl,GR,NTSh2014} we succeeded to obtain the following one (in an unessential
assumption $g^{\scriptscriptstyle00}=+1$)
\begin{eqnarray}
\lefteqn{
\hbox{\ssb E}\;=\;
\frac{\hbox{\ssb v}^{{{\prime}}{{\prime}}}\times(\hbox{\ssb s}
-s_{\scriptscriptstyle0}\hbox{\ssb v})}
{\bigl[(1+\hbox{\ssb v}{{\bcdot}}\hbox{\ssb v})(s_{\scriptscriptstyle0}{}^{2}
+\hbox{\ssb s}{{\bcdot}}\hbox{\ssb s})-(s_{\scriptscriptstyle0}+\hbox{\ssb s}{{\bcdot}}\hbox{\ssb v})^{2}\bigr]^{3/2}}
}\nonumber\\
&&\hfill\;-\;3\,\frac{(s_{\scriptscriptstyle0}{}^{2}+\hbox{\ssb s}{{\bcdot}}\hbox{\ssb s})
\;\hbox{\ssb v}^{{\prime}}\!{{\bcdot}}\hbox{\ssb v}-(s_{\scriptscriptstyle0}+\hbox{\ssb s}{{\bcdot}}\hbox{\ssb v})\;\hbox{\ssb s}{{\bcdot}}\hbox{\ssb v}^{{\prime}}  }
{\bigl[(1+\hbox{\ssb v}{{\bcdot}}\hbox{\ssb v})(s_{\scriptscriptstyle0}{}^{2}
+\hbox{\ssb s}{{\bcdot}}\hbox{\ssb s})-(s_{\scriptscriptstyle0}+\hbox{\ssb s}{{\bcdot}}\hbox{\ssb v})^{2}\bigr]^{5/2} }
\;\hbox{\ssb v}^{{\prime}}\times(\hbox{\ssb s}-s_{\scriptscriptstyle0}\hbox{\ssb v})
\hfill\nonumber\\
&&\hfill\label{35}\;+\;m\,\frac{(1+\hbox{\ssb v}{{\bcdot}}\hbox{\ssb v})\,\hbox{\ssb v}^{{\prime}}
-(\hbox{\ssb v}^{{\prime}}\!{{\bcdot}}\hbox{\ssb v})\,\hbox{\ssb v}}
{(1+\hbox{\ssb v}{{\bcdot}}\hbox{\ssb v})^{3/2}(s_{\scriptscriptstyle0}{}^{2}+\hbox{\ssb s}{{\bcdot}}\hbox{\ssb s})^{3/2}}\;.
\end{eqnarray}

Applying the prescription (\ref{29}) the above expression (\ref{35}) produces the corresponding
four-dimensional expression and thus the desired Euler-Poisson equations
\begin{equation}\label{36}
\fbox{$
\displaystyle
{\bmat{\mathcal E}}=
\frac{\ast\,\bmat{\ddot{ u}}\wedge\bmat u\wedge\bmat s}
{\|\bmat s\wedge\bmat u\|^3}
\;-\;3\,\frac{\ast\,\bmat{\dot{ u}}\wedge\bmat u\wedge\bmat s}
{\|\bmat s\wedge\bmat u\|^5}\,(\bmat{\dot{ u}}\wedge\bmat s){{\bcdot}}(\bmat u\wedge\bmat s)
+\frac m{\|{\bmat s}\|^3}\left[\frac{\bmat{\dot{ u}}}{\|\bmat u\|}
\,-\,\frac{\bmat{\dot{ u}}{{\bcdot}}\bmat u}{\|\bmat u\|^3}\,\bmat u
\right]=\boldsymbol0
$}
\end{equation}

Assume ${\x}_{_{\bmat 3}}$ denote the third-order prolongation of the generator
${\x}$ (\ref{22})
of (pseu\-do-\negthinspace\nolinebreak) \negthinspace\negthinspace or\-tho\-go\-nal transformations with the group parameters
$\hbox{\ssb n}$, $\hbox{\ssb q}$ to the space $J_{3}(\R;\R^3)$,
\begin{eqnarray}\label{37}
{\x}_{_{\bmat 3}}&=&-\,(\hbox{\ssb q}{{\bcdot}}\hbox{\ssb s})\,\partial_{_{{\scriptstyle s}^{\scriptscriptstyle0}}}
\,+\,s^{\scriptscriptstyle0}\,(\hbox{\ssb q}\,\bkey.\,\bpartial_{_{\eisf s}})
\,+\,\bkey[\hbox{\ssb n},\hbox{\ssb s},\bpartial_{_{\eisf s}}\bkey]
\,-\,(\hbox{\ssb q}{{\bcdot}}\hbox{\ssb x})\,\partial_{_{\scriptstyle t}}
\,+\,t\,(\hbox{\ssb q}\,\bkey.\,\bpartial_{_{\eisf x}})
\,+\,\bkey[\hbox{\ssb n},\hbox{\ssb x},\bpartial_{_{\eisf x}}\bkey]
\nonumber\\
&+&(\hbox{\ssb q}\,\bkey.\,\bpartial_{_{\eisf v}})
\,+\,(\hbox{\ssb q}{{\bcdot}}\hbox{\ssb v})\,(\hbox{\ssb v}\,\bkey.\,\bpartial_{_{\eisf v}})
\,+\,\bkey[\hbox{\ssb n},\hbox{\ssb v},\bpartial_{_{\eisf v}}\bkey]
\nonumber\\
&+&2\,(\hbox{\ssb q}{{\bcdot}}\hbox{\ssb v})\,(\hbox{\ssb v}^{{\prime}}\!\bkey.\,\bpartial_{_{\eisf v^\prime}})
\,+\,(\hbox{\ssb q}{{\bcdot}}\hbox{\ssb v}^{{\prime}})\,(\hbox{\ssb v}\,\bkey.\,\bpartial_{_{\eisf v^\prime}})
\,+\,\bkey[\hbox{\ssb n},\hbox{\ssb v}^{{\prime}}\!,\bpartial_{_{\eisf v^\prime}}\bkey]\,.
\end{eqnarray}

The following assertions are true:
{\em
\begin{enumerate}
\item
Vector quantity ${\bmat{\mathcal E}}$ in (\ref{36}) constitutes a system of Euler-Poisson expressions;
\item
Let vector field ${\x}_{_{\bmat 3}}$  be given by (\ref{37}), and let $\hbox{\ssb E}$
be given by (\ref{35}) (recall that (\ref{35}) is nothing more but merely (\ref{36}),
parametrized with respect to the coordinate time $t=x^{\scriptscriptstyle0}$).
The relation
$$
{\bf L}({\x}_{_{\bmat 3}})(\hbox{\ssb E})\;=\;\hbox{\ssb n}\times\hbox{\ssb E}
\;+\;(\hbox{\ssb q}{{\bcdot}}\hbox{\ssb v})\;\hbox{\ssb E}
\;-\;(\hbox{\ssb v}{{\bcdot}}\hbox{\ssb E})\;\hbox{\ssb q}
$$
proves the invariance of the Euler-Poisson equations
(\ref{36}), and at the same time shows that whatever a possible
Lagrange function for (\ref{36}) might exist, it by no means will reveal invariance
under (pseu\-do-\negthinspace\nolinebreak) \negthinspace\negthinspace or\-tho\-go\-nal transformations, even in the generalized sense (that
is up to a total derivative term).
\item
Equations (\ref{36}) describe (in metric signature $2$) the motion of relativistic spinning
free particle with constant spin four-vector $\bmat s$ and with the rest mass
$$
m_{\scriptscriptstyle0}\;=\;m\left[1\,-\,\frac{(\bmat s{{\bcdot}}\bmat u)^{2}}
{(\bmat s{{\bcdot}}\bmat s)\,(\bmat u{{\bcdot}}\bmat u)}\right]^{3/2}_{^{\ \ \ {\displaystyle.}}}
$$
\end{enumerate}
}

\subparagraph{A comment on the order of the Euler-Lagrange form.}
Euler-Lagrange equations are polynomial with respect to the derivatives of orders
greater than the maximal order of the derivatives which enter in the corresponding
Lagrange function. It was this property that inspired us to try to diminish
the order of the underlying manifold $J_{s}(\R,\R^{q})$ from $s=3$ to
$s=2$ by means of introducing  the differential form
${\underline{\bmat{\epsilon}}}$ in (\ref{33})
rather than ${\bmat{\epsilon}}$.

\end{document}